\documentclass[12pt]{gtart}
\usepackage{amscd,amssymb,frak,amsmath,pb-diagram}

\title{On Rigidly Scalar-Flat Manifolds} 

\author{Boris Botvinnik, Brett McInnes} 
\address{University of Oregon,  National University of Singapore} 
\email{botvinn@math.uoregon.edu, matmcinn@nus.edu.sg}   

\thanks{{\bf Date: November 3, 1999}}



\newcommand{\abstracttext}{Witten and Yau (hep-th/9910245) have
        recently considered a generalisation of the AdS/CFT
        correspondence, and have shown that the relevant manifolds
        have certain physically desirable properties when the scalar
        curvature of the boundary is positive.  It is natural to ask
        whether similar results hold when the scalar curvature is
        zero. With this motivation, we study compact scalar flat
        manifolds which do not accept a positive scalar curvature
        metric. We call these manifolds rigidly scalar-flat. We study
        this class of manifolds in terms of special holonomy
        groups. In particular, we prove that if, in addition, a
        rigidly scalar flat manifold $M$ is $Spin$ with $\dim M\geq
        5$, then $M$ either has a finite cyclic fundamental group, or
        it must be a counter example to Gromov-Lawson-Rosenberg
        conjecture.}

\newtheorem{Theorem}{Theorem}[section]

\newtheorem{Lemma}[Theorem]{Lemma}
\newtheorem{Corollary}[Theorem]{Corollary}

\theoremstyle{definition}

\newtheorem{Conjecture}[Theorem]{Conjecture}

\theoremstyle{remark}

\newtheorem{Remark}[Theorem]{Remark}
\def\Box{\qed}
\def\mathfrak{\frak}
\newcommand{\HK}{{\mathrm H}{\mathrm K}}
\newcommand{\vol}{{\mathrm v}{\mathrm o}{\mathrm l}}
\newcommand{\sym}{{\mathrm s}{\mathrm y}{\mathrm m}}
\newcommand{\Hol}{{\mathrm H}{\mathrm o}{\mathrm l}}
\newcommand{\Ker}{{\mathrm K}{\mathrm e}{\mathrm r}}
\newcommand{\Z}{{\mathbf Z}}

\newcommand{\R}{{\mathbf R}}
\newcommand{\CP}{{\mathbf C}{\mathbf P}}
\newcommand{\C}{{\mathbf C}}
\newcommand{\Cl}{{\mathcal C}\ell}
\begin{document} 
\begin{abstract} \abstracttext \end{abstract} 
\maketitlepage
\section{Introduction}\label{s1}
{\bf \ref{s1}/1. Motivation: Physics.} 
Witten and Yau \cite{WY} have
recently considered the generalization of the celebrated AdS/CFT
correspondence to the case of a complete $(n+1)$-dimensional Einstein
manifold $W$ of negative Ricci curvature with a compact conformal (in
the sense of Penrose) boundary $M$ of dimension $n$. (In the simplest
case, $W$ is a hyperbolic space, and $M$ is the standard sphere.) It
can be shown that the relevant conformal field theory on $M$ is stable
if the scalar curvature $R_M$ of $M$ is positive, but not if $R_M$ can
be negative. When $R_M=0$, stability is more delicate; the theory is
stable in some cases, but not in the others. Witten and Yau \cite{WY}
show that if $R_M$ is positive, then $W$ and $M$ have several
physically desirable properties: the boundary $M$ is connected, and
$W$ is free of wormholes. Furthermore, the topology of $M$ exerts a
strong influence on that of $W$ (the fundamental group of $W$ is no
``larger'' than that of $M$).

We wish to ask whether similar statements are valid when $R_M=0$ but
the conformal field theory on $M$ is still stable. First, suppose that
$M$ is such that the conformal structure can be perturbed so that the
scalar curvature becomes positive everywhere; it is obvious that the
Witten-Yau conclusions then hold. {\it Therefore, the only non-trivial
case is the one in which such perturbations are not possible.} We shall say
in the latter case that $M$ is {\it rigidly} scalar-flat. To summarize,
then, the {\it rigidly} scalar-flat boundaries are the only cases in which
the difficulties discussed by Witten and Yau could possibly arise.

For the benefit of physically oriented readers, we can expand on this
crucial point as follows. Let $M$ be any scalar-flat compact
manifold. Let $x$ be a point in $M$, and perturb the metric slightly
(but {\it not} conformally) so that the scalar curvature at $x$ is an
arbitrarily small positive number. This perturbation can give rise to
two possibilities: either it maintains $R_M\geq 0$ everywhere on $M$,
or it does not. (That is, in the second case, $M$ is such that the
scalar curvature {\it must} become negative at some other point if it
becomes positive at $x$.) In the latter case we say that $M$ is {\it
rigidly scalar flat}, while in the former it is non-rigid. Now suppose
that $M$ is {\it not} rigidly scalar-flat, and that we slightly
perturb as above so that $R_M\geq 0$ everywhere, and $R_M(x)>0$. Then
according to a theorem of Kazdan and Warner \cite{KW}, by means of a
conformal deformation of the metric we can force the scalar curvature
to be positive {\it everywhere}. (Recall that the metric on $M$ is in
any case only defined modulo conformal factors, so this procedure
changes nothing physical.) But now we apply the result of Witten and
Yau \cite{WY}. {\it We conclude that the results of Witten and Yau are
valid even when $R_M=0$, provided that $M$ is} not {\it rigidly
scalar-flat}.

At this point it is clear that we need two things: first, a way of
distinguishing compact rigidly scalar-flat manifolds from those which
are not rigid; and second, a way of extending the results of Witten
and Yau to the rigid case. 

The objective of this work is to supply the first of these. Our
results essentially state that, given the truth of the
Gromov-Lawson-Rosenberg conjecture, compact, locally irreducible,
rigidly scalar flat $Spin$ manifolds have structures which are very
special and can be described very precisely. We hope that this precise
description of these manifolds  will facilitate a completion of the
Witten-Yau results. 

The reader will notice in particular that the fundamental group of $M$
is severely constrained: it must be a cyclic group of low
order. Presumably this has important consequences for $W$. On a less
conjectural note, manifolds {\it not} in the lists given in Theorem
\ref{BR6} below are (given certain easily verified conditions) {\it
not} rigidly scalar flat. Thus, for example, the familiar
six-dimensional Calabi-Yau quintic in $\CP^4$ is obviously
scalar-flat, but it is not rigidly scalar flat: so the Witten-Yau
conclusions are valid in this case.

{\bf \ref{s1}/2. Motivation: Geometry.}  
Henceforth we
consider only smooth, closed and compact manifolds.  As above, a
scalar-flat Riemannian manifold is called {\it rigidly scalar-flat} if
it accepts no metric of positive scalar curvature.  The manifold is
``rigid'' in the sense that if a scalar-flat metric is perturbed so
that the scalar curvature becomes positive at some point, then it must
become negative at some other point. A basic theorem of Bourguignon
\cite[Corollary 4.49]{Besse} states that for a rigidly scalar-flat
manifold the space of Ricci-flat metrics coincides with the space of
scalar-flat metrics.

We recall a well-known fact: that all closed manifolds (of dimension at
least three) are divided into the following three classes (see
\cite{KW}):
\begin{enumerate}
\item[{\bf (P)}] manifolds accepting a positive scalar curvature metric;
\item[{\bf (Z)}] manifolds accepting a non-negative scalar curvature metric,
but not accepting a positive scalar curvature metric;
\item[{\bf (N)}] manifolds not accepting a metric of a non-negative
scalar curvature.
\end{enumerate}
It follows from \cite{KW} that the class of rigidly scalar-flat
manifolds is nothing but the class {\bf (Z)}.
Our objective in this paper is to study and emphasize general
properties of rigidly scalar-flat manifolds. 

Since these manifolds are necessarily Ricci-flat, it follows from the
Cheeger-Gromoll theorem that they have the structure $ (T^r
\times \widehat{M})/F, $ where $T^r$ is a flat torus, and $F$ is a
finite group and so is the fundamental group of $\widehat{M}$. By the
work of Gromov-Lawson \cite{GL2}, $T^r \times \widehat{M}$ is rigidly
scalar-flat if and only if $\widehat{M}$ is rigidly scalar-flat, so we
see that the study of the structure of compact rigidly scalar-flat
manifolds reduces to understanding the behavior of the rigidity condition
under {\it finite} coverings and quotients. It is therefore natural to
begin by removing the toral factor.

We shall consider manifolds which are compact and {\it locally
irreducible}. That is, the universal Riemannian cover is not isometric
to a product of lower-dimensional manifolds. If $\dim M> 1$, this
eliminates the toral factor, and forces the fundamental group
$\pi_1(M)$ to be finite. 
In addition, we assume that all manifolds we consider are oriented.

When the manifold is simply connected, and $\dim M\geq 5$, there are two
very different cases:  
\begin{enumerate}
\item[$\bullet$] $M$ is not a $Spin$ manifold,
\item[$\bullet$] $M$  is a $Spin$ manifold.
\end{enumerate}
In the first case the work of Gromov-Lawson \cite{GL} implies that
there are no rigidly scalar-flat manifolds, since a simply
connected non-$Spin$ manifold (with $\dim M\geq 5$) always accepts a
metric of positive scalar curvature.
In the $Spin$ case Stolz' theorem \cite{St1} says that $M$ (again, if
$\dim M\geq 5$) does not accept a positive scalar curvature metric if
and only if $\alpha(M)\neq 0$ in the real $K$-theory $KO_n$. Here,
recall that $\alpha$ is the index of the Dirac operator valued in the real
$K$-theory. In dimensions equal to a multiple of $4$ we use the
notation $\hat{A}$ instead of $\alpha$.

According to Besse, an $n$-dimensional locally irreducible Riemannian
manifold has {\it generic} holonomy if its restricted holonomy group is
isomorphic to $SO(n)$, or to $U({n\over 2})$, if it is locally
K\"ahlerian. Otherwise it has {\it special} holonomy. 
It was observed by Futaki \cite{Futaki} that
Stolz' theorem implies that a rigidly scalar-flat compact $Spin$
simply connected manifold has special holonomy. But Besse emphasizes
that such manifolds are Einstein, and, of course, a scalar-flat
Einstein manifold is Ricci-flat.
Thus we see that Bourguignon's theorem should be understood in terms
of the concept of special holonomy.
We believe that this phenomenon underlies Bourguignon's theorem
in general: that is, locally irreducible rigidly scalar-flat compact
manifolds are Ricci-flat {\bf because} they have special holonomy.
That is, we conjecture
\begin{Conjecture}\label{WeakBesse1}{\rm (Weak Besse Conjecture)}
Every locally irreducible rigidly scalar-flat compact manifold
has special holonomy.
\end{Conjecture}
\begin{Remark}
Besse mentions \cite[p. 19]{Besse} the possibility that
{\it all} locally irreducible compact Ricci-flat manifolds have
special holonomy. We call this statement Besse's Conjecture.
\end{Remark}
\section{Some basic results}
We emphasize that we consider here $Spin$ manifolds.  Now as a guide to
deciding what one can reasonably expect to prove, we have
\begin{Lemma}\label{BR1}
Let $M$ be a compact locally irreducible rigidly scalar-flat $Spin$
manifold, of dimension $\geq 5$. Suppose either that \\
\begin{tabular}{ll}
& {\bf (a)} $\dim M$ is not a multiple of $4$,\\
{\bf or} & 
{\bf (b)} $\dim M$ is a multiple of $4$, but $\pi_1(M)$ is not cyclic.
\end{tabular}\\
Then $M$ is a counter-example to the Gromov-Lawson-Rosenberg Conjecture.
\end{Lemma}
Before proving Lemma \ref{BR1}, we briefly recall a few things. Let
$\pi_1 M= \pi$ be a finite group; then a classifying map $f: M
\longrightarrow B\pi$ defines the element $[(M,f)]\in
\Omega^{Spin}_n(B\pi)$, where $\Omega^{Spin}_n(\cdot)$ is the $Spin$
cobordism theory. Then $\alpha$ is the composition
$$
\alpha: \Omega^{Spin}_n(B\pi) \stackrel{D}{\longrightarrow} KO_n (B\pi) 
\stackrel{A}{\longrightarrow} KO_n (C^*_r\pi),
$$
where $D$ is the Atiyah-Bott-Shapiro map, and $A$ is the assembly map;
see \cite{R,RS,RS2} for details. 

The $K$-theory groups $KO_* (C^*_r\pi)$ and the assembly map $A$
are well-understood when $\pi$ is a finite group;
see \cite{RS2}. In particular, the group $KO_n (C^*_r\pi)$ and the
invariant $\alpha(M)$ may be described as follows.  Since $\pi=\pi_1
M$ is a finite group, we consider the universal cover $\widetilde{M}
\longrightarrow M$.  A $Spin$ structure on $M$ (let $\dim M =n$)
lifts to $\widetilde{M}$, and gives the Dirac operator $D$ on
$\widetilde{M}$, which commutes with the action of the Clifford algebra
$\Cl_n$ and the deck-transformations of $\pi$.  This turns the kernel
$\Ker(D)$ of the Dirac operator $D$ into a module over the algebra
$\Cl_n\otimes \R\pi$, where $\R\pi$ is the real group ring of
$\pi$. Let $i: \Cl_n \rightarrow \Cl_{n+1}$ be the natural inclusion,
and let ${\mathfrak G}(\Cl_j\otimes \R\pi)$ be the corresponding
Grothendieck group. Then the invariant $\alpha(M)$ is nothing but the
residue class of $\Ker (D)$ in the corresponding $K$-theory:
$$
\alpha(M) = [\Ker (D)]\in {\mathfrak
G}(\Cl_n \otimes \R\pi)/i^*{\mathfrak G}(\Cl_{n+1} \otimes \R\pi)= 
KO_n(\R\pi)= KO_n (C^*_r\pi).
$$
Thus the $K$-theory $KO_n(\R\pi)$ is ``assembled'' out of regular
$KO_*$-groups, parametrized by representations of $\pi$. In
particular, if $\alpha(M)\neq 0$ it means that there exists a harmonic
spinor on $\widetilde{M}$. Taking into account that the
Gromov-Lawson-Rosenberg Conjecture fails for non-finite fundamental
groups, we state it in the following form:
\begin{Conjecture}(Gromov-Lawson-Rosenberg Conjecture) A compact $Spin$
manifold $M$ with $\dim M \geq 5$ and finite fundamental group
does not accept a positive scalar curvature metric if and only if
$\alpha(M)\neq 0$.
\end{Conjecture}
\begin{Remark} Under the conditions of Lemma \ref{BR1} a manifold $M$
is compact, locally irreducible, and Ricci-flat,so $\pi_1(M)$ is
finite. Thus the Gromov-Lawson-Rosenberg Conjecture predicts under
these circumstances that rigidity implies $\alpha(M)\neq 0$.
\end{Remark}
\noindent
{\bf Proof of Lemma \ref{BR1}. Case (a)} Assume $M$ is not a
counter-example. Then $\alpha(M)\neq 0$, and $M$ admits a parallel
spinor, and so it follows from Futaki's paper \cite[Proposition
2.3]{Futaki} that the holonomy must be special. Now a compact $Spin$
manifold of dimension $n\neq$ multiple of $4$ always has $\alpha(M)=0$
except in dimensions $8k+1$ and $8k+2$. There are no scalar-flat
locally irreducible manifolds of special holonomy in odd dimensions
other than $7$ (by the Berger-Simon theorem (see \cite[Chapter
10]{Besse}) and the fact that scalar-flat symmetric spaces are
flat). Finally, it follows from \cite{Brett-5} that a (compact,
orientable) scalar-flat locally irreducible  manifold of special
holonomy in dimension $8k+2$ is necessarily of holonomy $SU(4k+1)$.

Note that this is the {\it full} holonomy group, not just its identity
component, even if $M$ is not simply connected. The point here is that
this is a special property of dimensions $\neq $ multiple of $4$. If
we say that the {\it restricted} holonomy group of a manifold is
$SU(2)$ or $SU(4)$, this {\it does not mean} that it is a complex
manifold. But in dimension $10$, if the {\it restricted} holonomy is
$SU(5)$, then the manifold {\it has to be a K\"ahler manifold}. The
paper \cite{Brett-3} discusses this issue in detail.

Such a manifold, whether it is simply connected or not, is a complex
manifold (see \cite{W1,W2}), and by \cite[Theorem 1]{W1} and the main
result of \cite{W2}, the space of harmonic spinors (which coincides
with the space of parallel spinors in our case since the manifolds are
Ricci-flat), has complex dimension $2$. Recall (see \cite[Chapter III,
\S 10]{Lawson} that in dimension $2$ the index of the Dirac
operator $D$ is nothing but $\dim_{\C} (\Ker D)$ (mod 2) since 
$\Cl_1\cong \C$. Bott periodicity gives the same fact for all
dimensions $8k+2$. This would complete the proof in the case when
the original manifold $M$ is simply connected, so that $M =
\widetilde{M}$. However, this case is very special due to
\cite[Theorem 1]{W1}. Indeed, the dimension of the space of harmonic
spinors is equal to $2$ even after factorization by a finite group.
In particular, it means that for any finite group $\pi$ acting
holomorphically on a K\"ahler manifold $\widetilde{M}$ (the universal
cover of $M$), the kernel $\Ker D$ is always a trivial
$\R\pi$-module. Thus $\alpha(M)=0$.
So we have a contradiction in Case (a).

{\bf Case (b).} If $\alpha(M)\neq 0$ then, again, the holonomy group is
special. But a Ricci-flat compact locally irreducible $Spin$ manifold
of special holonomy with $n=$ a multiple of $4$ {\bf must} have a cyclic
fundamental group by the theorem of B. McInnes \cite{Brett-1}. Again,
we have a contradiction. \hfill $\Box$
%
\begin{Theorem}\label{BR2}
Let $M$ be a compact locally irreducible rigidly scalar-flat $Spin$
manifold with $\dim M \geq 5$. Then $M$ belongs to {\bf exactly} one
of the following classes.
\begin{enumerate}
\item[{\bf (I)}] $\alpha(M)=0$, and yet $M$ has a finite fundamental
group and admits no metric of positive scalar curvature. That is, $M$
violates the Gromov-Lawson-Rosenberg Conjecture.
\item[{\bf (II)}] $M$ has $\dim M = 4k$, and $\pi_1(M)$ is finite cyclic.
\end{enumerate}
\end{Theorem}
\noindent
{\bf Proof.} In view of Lemma \ref{BR1}, we only need to prove
that the two classes are disjoint. However, it is known \cite{BGS}
that the Gromov-Lawson-Rosenberg Conjecture is valid when $\pi_1(M)$
is finite cyclic. Hence no rigidly scalar-flat $Spin$ manifold can belong to
both classes. \hfill $\Box$
\begin{Corollary}\label{BR3}
The Weak Besse Conjecture for $Spin$ manifolds is valid if $\dim M =
4k>4$ and $\pi_1(M)$ is finite cyclic.
\end{Corollary}
\begin{Corollary}\label{BR4}
Let $\pi$ be a finite non-cyclic group such that it can be shown that
for $Spin$ manifolds of fundamental group $\pi$, there exists a metric
of positive scalar curvature if and only if $\alpha(M)=0$. Let $M$ be
a compact locally irreducible $Spin$ manifold of dimension $\geq 5$
with $\pi_1(M)=\pi$. Then $M$ cannot be rigidly scalar flat.
\end{Corollary}
\begin{Remark}
The Gromov-Lawson-Rosenberg Conjecture is known to be true when a
manifold has fundamental group belonging to a short list of groups,
see \cite{RS} for details.  For instance, this is the case when
fundamental group is a space form group, i.e. cyclic, quaternionic or
generalized quaternionic, see \cite{BGS}.
\end{Remark}
\begin{Corollary}\label{BR5}
Let \ $M$ \ be any compact locally irreducible \ $Spin$ \ manifold
with $\dim M \geq 5$, $\dim M \neq 4k$ and $\alpha(M)\neq 0$. Then the
scalar curvature is negative at some point in $M$.
\end{Corollary}
\begin{Remark}
In other words, these manifolds belong to the class {\bf
(N)}.  There are many examples of manifolds like this, for instance,
Hitchin's exotic spheres \cite{Hitchin}, or B\'erard-Bergery's
nine-dimensional example with $\pi_1(M)=\Z_2$; see \cite{RS2}.
\end{Remark}
Theorem \ref{BR2} means that compact locally irreducible
rigidly scalar-flat $Spin$ manifolds fall into two disjoint classes:
the exotic Class (I), and the Class (II). About the Class (I) we
can say little, except that we do not believe it is nonempty. About Class
(II), by contrast, we can say a great deal. We introduce some notation.
Let $n$ be a multiple of $4$. Then we define the following finite set of 
odd numbers:
$$
{\mathcal R}(n) =\left\{\begin{array}{cl}
\{\mbox{odd divisors of $2k+2$}\} & \mbox{if $n=8k+4$}
\\
\{\mbox{divisors of $2k+1$}\} & \mbox{if $n=8k$}.
\end{array}\right.
$$
We will denote by $\C_{\vol}$ a compact simply connected K\"ahlerian complex
manifold with a complex volume form $\omega$, and no other non-zero
holomorphic form (apart from multiples of $\omega$), and by
$\C_{\sym}$ a compact simply connected K\"ahlerian complex manifold
with a unique (up to constant scalar multiples) complex symplectic
form.
\begin{Theorem}\label{BR6}
The following classes of manifolds are identical.
\begin{enumerate}
\item[{\bf (1)}] Compact locally irreducible rigidly scalar-flat
	$Spin$ manifolds $M$ with $\dim M =n\geq 5$ equal to a multiple of $4$
	and with $\pi_1(M)$ cyclic.
\item[{\bf (2)}] Compact $n$-dimensional manifolds, where $n\geq 5$ is a
	multiple of $4$, with linear holonomy from the following list:
\begin{enumerate}
\item[{\bf (i)}] $SU({n\over 2})$.
\item[{\bf (ii)}] $SU({n\over 2}) \rtimes \Z_2$, if $n$ is a multiple
	of $8$, where the generator of $\Z_2$ acts by complex conjugation.
\item[{\bf (iii)}] $\Z_r \times Sp({n\over 4})$ with $r\in {\mathcal R}(n)$.
\item[{\bf (iv)}] $Spin(7)$, if $n=8$.
\end{enumerate}
\item[{\bf (3)}] Compact $n$-dimensional manifolds with $n\geq 5$ a
	multiple of $4$, with structure
\begin{enumerate}
\item[{\bf (i)}] $\C_{\vol}$.
\item[{\bf (ii)}] $\C_{\vol}/\Z_2$, $n=8k$, $\Z_2$ is generated by an
	antiholomorphic map.
\item[{\bf (iii)}] $\C_{\sym}/\Z_r$ with $r\in {\mathcal R}(n)$,
	$\Z_r$ is generated by a holomorphic map.
\item[{\bf (iv)}] $M$ is an $8$-dimensional simply connected manifold
	with a closed admissible $4$-form (see {\rm \cite{Joyce}}) and
	such that its Betti numbers satisfy $b^3 + b^4_+ = b^2 + b^4_-
	+ 25$.
\end{enumerate}
\end{enumerate}
\end{Theorem}
\section{Proofs and Remarks}\label{s2}
\noindent
{\bf Proof of Theorem \ref{BR6}.}
{\bf (1) $\Longrightarrow$ (3).} By Theorem \ref{BR2}, $\alpha(M)\neq
0$, so there is a parallel spinor, which means that the holonomy is
special, which implies that $M$ is equipped with an Einstein metric,
which, finally, must be Ricci flat. By \cite[Theorem 2]{Brett-1} $M$
is in one of the following classes:
\begin{enumerate}
\item[{\bf (i)}] $M$ is a Calabi-Yau manifold, that is, has holonomy
	$SU({n\over 2})$. When $n$ is a multiple of $4$, such a
	manifold must be simply connected (see \cite[p.113]{Salamon}),
	complex and Ricci-flat. By Bochner theory, any holomorphic form
	is parallel, hence $SU({n\over 2})$-invariant. By the
	representation theory of $SU({n\over 2})$, this means that
	there can be no non-zero holomorphic form other than a complex
	volume form and its constant multiples. As $M$ is simply
	connected and Ricci flat, the canonical bundle is trivial, so
	a complex volume form does indeed exist. Thus $M$ has
	structure $\C_{\vol}$.
\item[{\bf (ii)}] $M$ is a Calabi-Yau manifold factored by $\Z_2$,
	generated by an antiholomorphic map. Hence it has structure
	$\C_{\vol}/\Z_2$. This case occurs only if $n$ is a multiple
	of $8$.
\item[{\bf (iii)}] $M$ is $M^{\HK}/\Z_r$, $r\in {\mathcal R}(n)$,
	where $M^{\HK}$ is hyperK\"ahler. It is well known that
	compact hyperK\"ahler manifolds are simply connected complex
	manifolds (see \cite{Bea}). As the holonomy group of $M^{\HK}$
	is $Sp({n\over 4})$, and as any holomorphic form must be
	$Sp({n\over 4})$-invariant, the complex symplectic form must
	be unique up to constant multiples. Finally,
\item[{\bf (iv)}] $M$ could be $8$-dimensional manifold with holonomy
	$Spin(7)$. Joyce \cite{Joyce} shows that the Betti numbers
	satisfy the given relation.
\end{enumerate}

\noindent
{\bf (3) $\Longrightarrow$ (2).} 
\begin{enumerate}
\item[{\bf (i)}] Let $M$ be a compact simply connected complex
	manifold with no non-zero holomorphic form other than a
	complex volume form $\omega$.  
	Since $M$ is K\"ahlerian and 
	$\omega$ trivializes the canonical bundle, Yau's theorem
	\cite{Besse} gives us a Ricci-flat K\"ahler metric on $M$. As
	$\omega$ is the only non-zero holomorphic form, $M$ is
	irreducible, and so its holonomy is $SU({n\over 2})$.
\item[{\bf (ii)}] Let $M=\C_{\vol}/\Z_2$. As above, $\C_{\vol}$ has
	holonomy $SU({n\over 2})$, but since the generator of $\Z_2$
	does not act holomorphically, $M$ is not a K\"ahler manifold,
	and so its holonomy group is not $SU({n\over 2})$. As
	$\C_{\vol}$ is simply connected, $\pi_1(M)=\Z_2$ so the
	holonomy group has precisely $2$ components. A straightforward
	exercise in Lie theory shows that the only two-component
	subgroup of $SO(n)$ which has $SU({n\over 2})$ as identity
	component (and which is not contained in $U({n\over 2})$) is
	$SU({n\over 2})\rtimes \Z_2$, where the generator of $\Z_2$
	acts by complex conjugation.
\item[{\bf (iii)}] If the manifold is K\"ahlerian, and if $\sigma$ is
	the complex symplectic form, $\sigma^{{n\over 4}}$ trivializes
	the canonical bundle, so by Yau's theorem $\C_{\sym}$ admits a
	Ricci-flat K\"ahler metric. Since $\sigma$ is unique up to
	multiples, the holonomy group is $Sp({n\over 4})$. By an
	averaging argument combined with the Calabi uniqueness theorem
	(see \cite{Brett-4,Brett-5} for details) one can assume that
	the metric projects to $\C_{\sym}/\Z_r$, which becomes a
	K\"ahler manifold with $Sp({n\over 4})$ as the restricted
	holonomy group. Thus the holonomy group is $(\Z_q \times
	Sp({n\over 4}))/ \Z_2$ for some $q$. Since a $\C_{\sym}$
	manifold is simply connected, $\pi_1(\C_{\sym}/\Z_r)=\Z_r$, so
	there is a homomorphism from $\Z_r$ onto $\Z_{{q/2}}$ (we can
	assume that $q$ is even without loss of generality). Thus
	${q/2}$ divides $r$, which is odd, so ${q/2}$ is odd; hence
	the holonomy group is actually $\Z_{{q/2}}\times Sp({n\over
	4})$. Now suppose ${q/2}<r$, and let $k=r/({q\over 2})$.  Then
	$\C_{\sym}/\Z_k$ is a ${q\over 2}$-fold covering of
	$\C_{\sym}/\Z_r$, and so its full holonomy group is precisely
	$Sp({n\over 4})$. Now let $f : \C_{\sym} \longrightarrow
	\C_{\sym}$ generate $\Z_r$, so $f^{{q\over 2}}$ is a
	non-trivial fixed-point-free holomorphic map generating
	$\Z_k$. Since $\Hol(\C_{\sym}/\Z_k)=Sp({n\over 4})$,
	$f^{{q\over 2}}$ preserves $\sigma$ and all of its powers. But by
	the representation theory of $Sp({n\over 4})$ (see \cite{Bea})
	$\sigma$ and its powers are the {\it only} non-zero
	holomorphic forms on $\C_{\sym}$. Since all these forms are of
	{\it even} degree, the holomorphic Lefschetz number of
	$f^{{q\over 2}}$ is non-zero, which is impossible. Thus
	actually ${q\over 2}=r$ and we have that
	$\Hol(\C_{\sym}/\Z_r)=\Z_r\times Sp({n\over 4})$.
\item[{\bf (iv)}] Let $M$ be as described. Then Joyce \cite{Joyce}
	shows that $M$ has a metric of holonomy $Spin(7)$.
\end{enumerate}
{\bf (2) $\Longrightarrow$ (1).}
\begin{enumerate}
\item[{\bf (i)}] Let $M$ be a compact $n$-dimensional manifold of
	holonomy $SU({n\over 2})$, where $n\neq 4$ is a multiple of
	$4$. Then $M$ is scalar-flat, simply connected, and $Spin$. It
	is well-known that such manifolds are rigid (since
	$\hat{A}(M)=2\neq 0$). 
\item[{\bf (ii)}] Manifolds $M$ of holonomy $SU({n\over 2})\rtimes
	\Z_2$ are double-covered by manifolds of holonomy $SU({n\over
	2})$, so it is evident from (i) that they must be scalar-flat,
	rigid, and have $\pi_1(M)=\Z_2$, which is indeed cyclic. It
	remains only to prove that these manifolds are $Spin$. For
	this, observe that {\it when $n$ is a multiple of $4$}, the
	preimage of $SU({n\over 2})\rtimes \Z_2$ in $Spin(n)$ is
	$\{\pm 1\}\times SU({n\over 2})\rtimes \Z_2$, so $SU({n\over
	2})\rtimes \Z_2$ is a subgroup of $Spin(n)$. Thus if $H(M)$ is
	a bundle of orthonormal frames with $SU({n\over 2})\rtimes
	\Z_2$ as structural group, we can define
$$
	S(M) = {H(M)\times SO(n)\over SU({n\over 2})\rtimes \Z_2},
\ \ \ \mbox{with} \ \ \ S(M)/\{\pm 1\} = SO(M),
$$
	the corresponding {\it full} bundle of orthonormal frames. So
	$M$ is $Spin$. (Actually, $\hat{A}(M)=1$. See \cite{MS} for
	details of this argument.)
\item[{\bf (iii)}] Manifolds of holonomy $\Z_r\times Sp({n\over 4})$
	are of the form $M = \hat{M}/F$ for some manifold $\hat{M}$ of
	holonomy $Sp({n\over 4})$. It is well known that $\hat{M}$ is
	scalar-flat and rigid, hence the same is true of $M =
	\hat{M}/F$. Clearly $F$ acts freely and holomorphically on the
	hyperK\"ahler manifold $\hat{M}$, but (as is shown in
	\cite{Brett-6}) this implies that $F$ is cyclic. Finally, the
	preimage of $\Z_r\times Sp({n\over 4})$ in $Spin(n)$ is $\{\pm
	1\}\times \Z_r \times Sp({n\over 4})$, so, as above, these
	manifolds are $Spin$. (Note that $\hat{A}(M)=(1+{n\over
	4})/r$.)
\item[{\bf (iv)}] Compact $8$-dimensional manifolds with holonomy
	$Spin(7)$ are necessarily simply connected, $Spin$, and
	rigidly scalar-flat (since $\hat{A}=1$).
\end{enumerate}
This completes the proof. \hfill $\Box$
\begin{Remark}
	The $Spin$ condition plays an essential role throughout this subject.  For
	example, in four dimensions, there are compact manifolds which
	are not $Spin$ and yet are rigidly scalar-flat (e.g. any Enriques
	surface). Similarly there are rigidly scalar-flat manifolds
	with non-cyclic fundamental groups (e.g the Hitchin manifold
	$K3/(\Z_2\times \Z_2)$.)  Such examples are known in all
	dimensions equal to a multiple of $4$; see \cite{Brett-4} for their
	explicit description.
\end{Remark}
\begin{Remark} When $n=8$ the possibilities for holonomy groups are
\begin{equation}\label{hol-8}
\begin{array}{ccccl}
SU(4) & = & Spin(6)& \subset& Spin(7),
\\
SU(4)\rtimes \Z_2 & = & Spin(6)\sqcup Spin(6)\cdot e_4e_5e_6e_7
&\subset& Spin(7),
\\
Sp(2)& = & Spin(5) &\subset& Spin(7),
\\
\Z_3\times Sp(2) &\subset& Spin(2)\cdot Spin(5) &\subset& Spin(7),
\\
Spin(7). &&&&
\end{array}
\end{equation}
At present there are known examples of manifolds with all of these holonomy
groups (\ref{hol-8}) except one: $\Z_3\times Sp(2)\subset Spin(7)$. 
\end{Remark}
We suspect that the techniques developed by Joyce \cite{Joyce} can be
extended to give a solution of the following
\smallskip

\noindent
{\bf Open Problem.} Construct a closed eight-dimensional manifold with 
holonomy group $\Z_3\times Sp(2)$. 
\begin{Remark} 
In the four-dimensional case, the situation is quite different, but
again there are obstructions to the existence of metrics of positive
scalar curvature; as is well known, these are associated with
Seiberg-Witten theory. The following result (for details, see \cite{T,
LeB, RS3}) is relevant.
\end{Remark}
\begin{Theorem}{\rm (Taubes \cite{T},  LeBrun \cite{LeB})}
Let $M$ be a closed, connected, oriented four-manifold with
$b_2^+(M)>1$. If $M$ admits a symplectic structure (in particular, if
$M$ is a K\"ahler manifold), then $M$ does not accept a positive
scalar curvature metric.
\end{Theorem}



\begin{thebibliography}{30}
\bibitem{Bea} A. Beauville, Vari\'et\'es K\"ahleriennes dont la
	premi\`ere classe de Chern est nulle. J. Diff. Geom. 18
	(1983), no. 4, 755--782
\bibitem{Besse} A. Besse, Einstein manifolds, Ergebnisse der
	Mathematik und ihrer Grenzgebiete (3) [Results in Mathematics
	and Related Areas (3)], 10. Springer-Verlag, Berlin-New York,
	1987.
\bibitem{BGS} B. Botvinnik, P. Gilkey, S. Stolz, The
	Gromov-Lawson-Rosenberg conjecture for groups with periodic
	cohomology. J. Diff. Geom. 46 (1997), no. 3, 374--405.
\bibitem{Futaki} A. Futaki, Scalar-flat closed manifolds not admitting
	positive scalar curvature metrics. Invent. Math. 112 (1993),
	no. 1, 23--29.
\bibitem{GL} M. Gromov, H. B. Lawson, The classification of simply
	connected manifolds of positive scalar curvature. Ann. of
	Math. (2) 111 (1980), no. 3, 423--434.
\bibitem{GL2} M. Gromov, H. B. Lawson, Spin and scalar curvature in
	the presence of a fundamental group. I. Ann. of Math. (2) 111
	(1980), no. 2, 209--230.
\bibitem{Hitchin} N. Hitchin, Harmonic spinors. Adv. in Math. 14
	(1974), 1--55.
\bibitem{Joyce} D. Joyce, Compact $8$-manifolds with holonomy ${\rm
	Spin}(7)$.  Invent. Math. 123 (1996), no. 3, 507--552
\bibitem{KW} J. L. Kazdan, F. W. Warner, Scalar curvature and
	conformal deformation of Riemannian structure. J. Differential
	Geometry 10 (1975), 113--13
\bibitem{Lawson} H. B. Lawson, M-L. Michelsohn, Spin geometry.
	Princeton Mathematical Series, 38.  Princeton University
	Press, Princeton, NJ, 1989. xii+427 pp.
\bibitem{LeB} C. LeBrun, On the scalar curvature of complex
	surfaces. Geom. Funct. Anal. 5 (1995), no. 3, 619--628.
\bibitem{Brett-1} B. McInnes, $Spin$ holonomy of Einstein
	manifolds. Comm. Math. Phys.  203 (1999), no. 2, 349--364.
\bibitem{Brett-2} B. McInnes, Metric symmetries and spin asymmetries
	of Ricci-flat Riemannian manifolds. J. Math. Phys. 40 (1999),
	no. 3, 1255--1267.
\bibitem{Brett-3} B. McInnes, Methods of holonomy theory for
	Ricci-flat Riemannian manifolds. J. Math. Phys. 32 (1991),
	no. 4, 888-896.
\bibitem{Brett-4} B. McInnes, Examples of Einstein manifolds with all
	possible holonomy groups in dimensions less than
	seven. J. Math. Phys. 34 (1993), no. 9, 4287--4304. 
\bibitem{Brett-5} B. McInnes, Holonomy groups of compact Riemannian
	manifolds: a classification in dimensions up to
	ten. J. Math. Phys. 34 (1993), no. 9, 4273--4286.
\bibitem{Brett-6} B.  McInnes, The quotient construction for a class
	of compact Einstein manifolds. Comm. Math. Phys. 154 (1993),
	no. 2, 307--312.
\bibitem{MS} A. Moroianu, U. Semmelmann. Parallel spinors and holonomy
       	groups. Preprint, Differential Geometry: math.DG/9903062.
\bibitem{R} J. Rosenberg, $C\sp *$-algebras, positive scalar curvature
	and the Novikov conjecture. II. Geometric methods in operator
	algebras (Kyoto, 1983), 341--374, Pitman Res. Notes
	Math. Ser., 123, Longman Sci. Tech., Harlow, 1986.
\bibitem{RS} J. Rosenberg, S. Stolz, Manifolds of positive scalar
	curvature.  Algebraic topology and its applications, 241--267,
	Math. Sci. Res. Inst. Publ., 27, Springer, New York, 1994.
\bibitem{RS2} J. Rosenberg, S. Stolz, A "stable" version of the
	Gromov-Lawson conjecture. The \v{C}ech centennial (Boston, MA,
	1993), 405--418, Contemp. Math., 181, Amer. Math. Soc.,
	Providence, RI, 1995.
\bibitem{RS3} J. Rosenberg, S. Stolz, Metrics of positive scalar
	curvature and connections with surgery, to appear,  
	http://www.math.umd.edu/$\sim$jmr/jmr$_{\mbox{-}}$pub.html 
\bibitem{Salamon} S. Salamon, Riemannian geometry and holonomy
	groups. Pitman Research Notes in Mathematics Series,
	201. Longman Scientific \& Technical, Harlow; copublished in
	the United States with John Wiley \& Sons, Inc., New York,
	1989.
\bibitem{St1} S. Stolz, Simply connected manifolds of positive scalar
	curvature, Ann.  of Math. (2) 136 (1992), no. 3, 511--540.
\bibitem{T} C. Taubes, The Seiberg-Witten invariants and symplectic
	forms. Math. Res. Lett. 1 (1994), no. 6, 809--822.
\bibitem{W1} M. Y. Wang, On non-simply connected manifolds with
	non-trivial parallel spinors. Ann. Global Anal. Geom.  13
	(1995), no. 1, 31--42.
\bibitem{W2} M. Y. Wang, Parallel spinors and parallel
	forms. Ann. Global Anal. Geom. 7 (1989), no. 1, 59--68.
\bibitem{WY} E. Witten, S.-T. Yau, Connectedness of the boundary in
	the AdS/CFT correspondence, Preprint, High-Energy Physics:
	hep-th/9910245.
\end{thebibliography}
\end{document}